# Statistical modeling for experiments with sliding levels

**Shao-Wei Cheng[1], C. F. J. Wu[2] and Longcheen Huwang[3]**

*Academia Sinica, Georgia Institute of Technology and National Tsing-Hua University*

**Abstract:** Design of experiment with related factors can be implemented by using the technique of sliding levels. Taguchi (1987) proposed an analysis strategy by re-centering and re-scaling the slid factors. Hamada and Wu (1995) showed via counter examples that in many cases the interactions cannot be completely eliminated by Taguchi's strategy. They proposed an alternative method in which the slid factors are modeled by nested effects. In this work we show the inadequacy of both methods when the objective is response prediction. We propose an analysis method based on a response surface model, and demonstrate its superiority for prediction. We also study the relationships between these three modeling strategies.

## 1. Introduction

In many investigations, the experimenters can choose an appropriate interval as the experimental range for each factor. The overall experimental region is then the cube formed by the tensor product of these intervals. Such an experimental region is called *regular*. However, when some of the factors are related, an appropriate experimental region becomes irregular and thus cannot be constructed in the usual manner. Factors are called *related* when the desirable experimental region of some factors depends on the level settings of other factors. Design of experiments with related factors can be implemented by using the technique of *sliding levels* proposed by Taguchi [7]. It has been used in practice for a long time but has received scant attention in the statistical literature. Some examples can be found in [2, 6, 7]. Li et al. [4] proposed a two-stage strategy for the sliding-level experiments whose desriable experimental region is unknown and needs to be explored during the experiment. Here the use of sliding is more complicated due to its engineering needs.

In this article we study the situations in which only one factor is chosen to be slid. This article is organized as follows. In Section 2, we will review the existing work on the sliding level technique and show the inadequacy of these methods when the objective of the experiment is response prediction. In Section 3, we will propose an analysis method based on a response surface model, and demonstrate its superiority for prediction. In Section 4, an illustration with a welding experiment will be given.

[1]Academia Sinica, Institute of Statistical Science, Taipei, Taiwan, e-mail: `swcheng@stat.sinica.edu.tw`
[2]Georgia Institute of Technology, School of Industrial and Systems Engineering, Atlanta, GA 30332-0205, USA, e-mail: `jeffwu@isye.gatech.edu`
[3]National Tsing-Hua University, Institute of Statistics, Hsinchu, Taiwan, e-mail: `huwang@stat.nthu.edu.tw`

*AMS 2000 subject classifications:* primary 62K15, 62K20; secondary 62P30.

*Keywords and phrases:* irregular experimental region, nested effect, re-centering and re-scaling, response prediction, response surface modeling, robust parameter design.

In Section 5, some results are presented based on a comparison between the response surface approach and Taguchi's approach. A summary is given in the last section.

## 2. Existing approaches

Taguchi [7] justified the use of sliding levels by the rationales of *bad region avoidance* and *interaction elimination*. The analysis strategy in his approach for sliding levels can be interpreted as a *re-centering* and *re-scaling* (RCRS) transformation, which transforms an irregular experimental region into a regular one as shown in Fig 1. In data analysis, this transformation is essentially to code the factor levels by regarding the slid factor as a non-slid factor. For example, $(+1, -1)$ is used for the conditional low and high levels respectively in a two-level slid factors, and $(-1, 0, +1)$ for the conditional low, median, and high levels respectively in a three-level slid factors with equally spaced levels. Consider two factors $A$ and $B$, in which there are several sliding levels for $B$ at each level of $A$. It is easy to show that an interaction in the original factor space is eliminated after RCRS only if the relationship between the mean response $E(y)$ and factors $A$ and $B$ satisfies the relationship:

$$(1) \qquad E(y) = g_1(x_A) + g_2\left[\frac{x_B - c_B(x_A)}{r_B(x_A)}\right],$$

where $g_1$ and $g_2$ are two arbitrary functions and the $c$'s and $r$'s represent the centering and scaling constants with those for factor $B$ depending on factor $A$. Furthermore, to eliminate the interaction between $A$ and $B$ for mean response satisfying (1), a proper choice of sliding levels based on $c$'s and $r$'s is required. As pointed out via a counter example by Hamada and Wu [3], inadequately locating the sliding levels will not remove the interaction. Similarly, an inadequate choice of scale will not eliminate the interaction neither.

One can infer that the sliding levels must be chosen properly in order to eliminate a potentially removable interaction. To achieve this, one has to know the exact relationship between the factors and the mean response $E(y)$. Because this relationship is not available, an experiment needs to be carried out. Therefore the advantage of interaction elimination by using sliding levels is questionable. Even though the related factors' interactions can be removed by proper centering and scaling, important information like robustness may be masked (see [3], for more details).

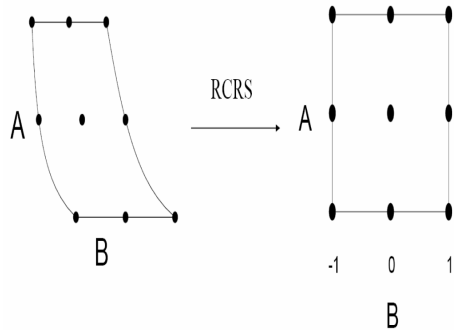

Fig 1. *Re-centering and re-scaling transformation of experimental region.*

Hamada and Wu [3] proposed a *nested-effects modeling* (NEM) approach by using a regression model with nested effects. Because the actual settings of the slid factor are different at each level combination of its related factors, sliding-levels designs can be viewed as nested designs. Hence, one can model the effect of the nested (slid) factor separately at each level combination of its related factors, i.e., the effects of the slid factor are defined conditional on the level combinations of its related factor. Consider the case of two related factors where factor $B$'s levels depend on $A$'s. The factor $A$ can be either qualitative or quantitative. For qualitative $A$, Hamada and Wu [3] proposed analyzing the effect of $B$ at each level of $A$. If $B$ is quantitative with more than two levels, the linear and quadratic effects of $B$ at the $i$th level of $A$ (denoted by $B_l|A_i$ and $B_q|A_i$) should be analyzed. Furthermore, the effects of factor $A$ are analyzed as well. For instance, if $A$ is qualitative with three levels, the two contrasts $A_{1,2}$ and $A_{1,3}$ can be considered, where $A_{i,j}$ represents the contrast between levels $i$ and $j$ of $A$, i.e., it denotes the difference between the average responses over the conditional levels of $B$ at level $i$ of $A$ and those at level $j$ of $A$. Because the levels of $B$ vary with the level of $A$, this is different from the usual meaning of $A_{i,j}$ in factorial designs with regular experimental region, where the same set of levels of $B$ is used for $i$ and $j$. If $A$ is quantitative, the linear and quadratic effects of $A$ (i.e., $A_l$ and $A_q$ in the linear-quadratic system defined in [8]) should be substituted for $A_{1,2}$ and $A_{1,3}$. The same reasoning will show that the meanings of $A_l$ and $A_q$ are again different from the usual ones.

The analysis using a regression model with nested effects resolves the problem that the sliding-levels design may not eliminate the interaction between related factors. It also provides more insight into the response-factor relationship and directly accounts for the relationship between related factors, which can be used to choose optimum factor levels. However, as far as response prediction is concerned, the nested effects analysis is incapable of accomplishing the task for quantitative $A$.

When $A$ is a quantitative factor, we may need to predict the response at a setting whose value of $A$, say $x_A^*$, is not included in the experimental plan. To achieve this, we need to have a fitted model of $B$ at $A = x_A^*$. However, such a model is not available in the NEM approach because an NEM offers fitted models of $B$ *only* for each levels of $A$ and $x_A^*$ is not one of the levels in the experiment. Therefore, response prediction at $x_A^*$ cannot be achieved in an NEM approach. Because the effects of $B$ are defined and analyzed conditional on $A$ in an NEM, $A$ is treated like a qualitative factor in the analysis about $B$. This results in the difficulty of response prediction at $x_A^*$. Turning to the RCRS approach for performing prediction, we have to know the centering and scaling constants of $B$ at $x_A^*$, i.e., $c_B(x_A^*)$ and $r_B(x_A^*)$, so that $B$ can be appropriately transformed at $A = x_A^*$ before substituting into the fitted RCRS model. However, both $c_B(x_A^*)$ and $r_B(x_A^*)$ may not be available to the investigators. In the next section, we shall propose an analysis method based on the response surface methodology and demonstrates its superiority for response prediction over the two existing approaches.

### 3. An analysis strategy based on response surface modeling

*Response surface modeling* (RSM) is an effective tool for building empirical models for the input and output variables in an experiment. In RSM, a true model is often expressed as $y = f(x_1, x_2, \ldots, x_k) + \epsilon$, where $y$ is the observed response, $f$ a function of $k$ quantitative factors $x_1, x_2, \ldots, x_k$, and $\epsilon$ an error term. For simplicity, the lowest level of a factor is coded as $-1$ and the highest level as $+1$. The function $f$ represents the response surface, which depicts the true relationship between the response and factors. Because the form of $f$ is often unknown, RSM replaces and approximates $f$ by a polynomial model of degree $d$ in the $x_i$'s. In practical applications, $d$ is often chosen to be one or two, and three when the response surface is expected to be more complicated and there are sufficient degrees of freedom. Fourth and higher degree polynomials are rarely used because they are not as effective as semi-parametric or nonparametric models. Further discussion on the response surface methodology can be found in [1] and [5].

In a sliding-level experiment, the adequate experimental region, denoted by $R_E$, usually has an irregular shape in contrast to the regular region in conventional factorial experiments. In such circumstances, the RSM can still be applied by first finding a cuboidal region that covers exactly the $R_E$ as follows. For each factor, let its lowest actual setting be coded as $-1$ and the highest actual setting as $+1$. Other settings of the factor is then proportionally coded according to their distances from the lowest one. In this coding, the cuboidal region $[-1, +1]^k$ is the smallest cube to cover the $R_E$. We call $[-1, +1]^k$ the *modeling region* and denote it by $R_M$. The RSM can then be applied in the modeling region to develop an empirical model. Unlike factorial designs with regular experimental region, the design points in a sliding-level experiment do not spread uniformly on the whole modeling region. Because there are no design points located in $R_M \setminus R_E$, we have no information about the response surface over $R_M \setminus R_E$. Therefore, the fitted model may fit well only in $R_E$, but not in the whole $R_M$. Another issue concerns the choice of appropriate polynomial models for the approximation of the true response surface. For sliding-level experiments, should we still use a $d$th-order polynomial model? This will be further explained later. When a fitted model is obtained, prediction can be easily done in the RSM approach. Its prediction is an interpolation in $R_E$ but an extrapolation in $R_M \setminus R_E$. An illustration of the RSM strategy will be given in Section 4.

Consider a nine-run experiment with factors $A$ and $B$, in which $A$ has three levels and conditional on each level of $A$, $B$ has three sliding levels. The NEM for the experiment can be written as:

$$(2) \qquad f(B|A_i) = b_0^i + b_l^i(B_l|A_i) + b_q^i(B_q|A_i), \quad i = -1, 0, 1,$$

where $b_0^i$, $b_l^i$ and $b_q^i$ are the conditional constant, linear, and quadratic main effects of $B$ given $A = i$. Because $A$ has three levels, the NEM has nine effects and therefore is saturated. On the other hand, a second-order RSM model for the experiment has only six effects. Because the NEM is saturated, it is clear that the RSM model is a submodel of the NEM. In other words, we can impose some constraints on the parameters of the NEM to obtain the RSM model. To find these constraints, we re-parameterize the NEM in (2) in terms of the coding based on the RSM as follows:

$$(3) \qquad f(x_B|x_A) = \alpha_{x_A} + \beta_{x_A} x_B + \gamma_{x_A} x_B^2, \quad x_A = -1, 0, 1,$$

where $x_B$ is coded according to the RSM approach but nested on $x_A$, and $\alpha_{x_A}$, $\beta_{x_A}$, and $\gamma_{x_A}$ are the zero-order, first-order, and second-order effects of $B$ conditional on $A = x_A$, respectively. Note that for $x_A = i$, $x_B$ is a linear transformation of $B_l|A_i$, and $x_B^2$ is a linear combination of 1, $B_l|A_i$, and $B_q|A_i$. By equating the NEM in (3) and the following second-order RSM model:

$$(4) \qquad f(x_A, x_B) \approx \lambda_0 + \lambda_1 x_A + \lambda_2 x_B + \lambda_{11} x_A^2 + \lambda_{22} x_B^2 + \lambda_{12} x_A x_B,$$

we obtain the following relationships:

(5)
$$\alpha_{x_A} = \lambda_0 + \lambda_1 x_A + \lambda_{11} x_A^2,$$
$$\beta_{x_A} = \lambda_2 + \lambda_{12} x_A,$$
$$\gamma_{x_A} = \lambda_{22}.$$

The equations in (5) indicate that the three conditional second-order effects of $B$ (i.e., $\gamma_i$'s) must be identical in the second-order RSM model, which save two degrees of freedom; the three conditional first-order effects of $B$ (i.e., $\beta_{x_A}$'s) must satisfy a linear constraint, which save one degree of freedom. The saving of three degrees of freedom explains why the RSM model has three parameters fewer than the NEM.

If the restrictions on $\beta_{x_A}$'s and $\gamma_{x_A}$'s in (5) are considered to be too rigid, we can add more parameters in the RSM model so that the corresponding $\beta_{x_A}$'s and $\gamma_{x_A}$'s can be free of the constraints as shown by the following relationships:

(6)
$$\alpha_{x_A} = \lambda_0 + \lambda_1 x_A + \lambda_{11} x_A^2,$$
$$\beta_{x_A} = \lambda_2 + \lambda_{12} x_A + \lambda_{112} x_A^2,$$
$$\gamma_{x_A} = \lambda_{22} + \lambda_{122} x_A + \lambda_{1122} x_A^2.$$

The resulting RSM model will be:

$$\begin{aligned} f(x_A, x_B) &\approx \lambda_0 + \lambda_1 x_A + \lambda_{11} x_A^2 \\ &\quad + (\lambda_2 + \lambda_{12} x_A + \lambda_{112} x_A^2) x_B \\ &\quad + (\lambda_{22} + \lambda_{122} x_A + \lambda_{1122} x_A^2) x_B^2 \\ &= \lambda_0 + \lambda_1 x_A + \lambda_2 x_B + \lambda_{12} x_A x_B + \lambda_{11} x_A^2 + \lambda_{22} x_B^2 \\ &\quad + \lambda_{112} x_A^2 x_B + \lambda_{122} x_A x_B^2 + \lambda_{1122} x_A^2 x_B^2. \end{aligned}$$

By adding three higher-order effects $x_A^2 x_B$, $x_A x_B^2$, and $x_A^2 x_B^2$ in the model (4), the RSM model has the same number of parameters and same capacity of estimation as the saturated NEM.

From the previous explanation, it is observed that the conventional RSM approach of using a $d$th-order model can be inappropriate for data from sliding-level experiments. For example, if a second-order model is adopted, some implicit constraints that can be impractical are placed on $\beta_{x_A}$'s and $\gamma_{x_A}$'s. However, if the experimenter would like to use a more complicated model, such as a third-order model, there are not enough degrees of freedom for estimating all parameters.

Another interesting observation about the relationship between NEM and RSM model can be obtained from the equations in (6). Consider, for example, the three conditional zero-order effects, $\alpha_{x_A}$'s. They are individually estimated at each level of $A$. From (6), $\alpha_{x_A}$'s can be expressed as a quadratic polynomial of $x_A$ with coefficients from parameters in the RSM model. To estimate these parameters, we can first estimate the $\alpha_{x_A}$'s, denoted by $\hat{\alpha}_{x_A}$, by least squares and then solve the equations $\hat{\alpha}_{x_A} = \lambda_0 + \lambda_1 x_A + \lambda_{11} x_A^2$, for $x_A = -1$, 0, 1, to obtain $\hat{\lambda}_0$, $\hat{\lambda}_1$, and $\hat{\lambda}_{11}$. In other words, for $x_A^*$ which is not in $\{-1, 0, 1\}$, we can predict $\alpha_{x_A^*}$ by using $\hat{\lambda}_0 + \hat{\lambda}_1 x_A^* + \hat{\lambda}_{11} (x_A^*)^2$. The same procedure can be applied to $\beta_{x_A}$'s and $\gamma_{x_A}$'s in (6).

It is then clear why and how the RSM model can be used for prediction. Suppose that we want to predict the value of $E(y)$ at $(x_A^*, x_B^*)$, where $x_A^*$ is not included in the experimental plan. From the argument given in Section 2, the NEM cannot be used for prediction at $x_A^*$ because no data are collected at $x_A^*$ for estimating the conditional effects $\alpha_{x_A^*}$, $\beta_{x_A^*}$, and $\gamma_{x_A^*}$. However, the RSM model treats $\alpha_{x_A}$, $\beta_{x_A}$, and $\gamma_{x_A}$ as continuous (second-order) polynomials over $x_A$. From this viewpoint and (6), the predicted value of $E(y)$ at $(x_A^*, x_B^*)$ is simply $\hat{\alpha}_{x_A^*} + \hat{\beta}_{x_A^*} x_B^* + \hat{\gamma}_{x_A^*} x_B^{*2}$, where $\hat{\alpha}_{x_A^*}$, $\hat{\beta}_{x_A^*}$, and $\hat{\gamma}_{x_A^*}$ are obtained by substituting $x_A^*$ into the right hand side expressions in (6) with $\lambda$'s replaced by $\hat{\lambda}$'s. Note that in the prediction procedure using RSM approach, the $\alpha_{x_A}, \beta_{x_A}$, and $\gamma_{x_A}$ are assumed to be *continuous* functions over $x_A$ and their changes over $x_A$ are assumed to follow the quadratic polynomials in (6). These assumptions explain why prediction is feasible in the RSM, but not in the NEM. In other words, the RSM approach regards the three levels of $A$ as quantitative and utilizes some continuity assumptions on $A$ for prediction. When similar assumptions are imposed on an NEM, prediction using NEM can be feasible.

We will show in Sections 4 and 5 that the RSM model for sliding-level experiment can suffer from severe collinearity between the effects of the slid factor and the effects of its related factors. The RSM model is therefore not a good choice for the purpose of identifying important effects, especially when it is required to perform model selection, such as forward selection or $C_p$. In these circumstances, we can adopt the following *hybrid strategy* that combines NEM and RSM as follows.

(i) It starts from a NEM, which has better orthogonality between effects in the models.
(ii) After important effects are identified, we can translate the fitted NEM into an RSM model through equations that relate the parameters in the two models (such as (6)). The resulting RSM model can then be used for response prediction.

## 4. Illustration: a welding experiment

We illustrate the three modeling strategies and compare their results by using data from a welding experiment reported in Chen, Ciscon, and Ratkus [2]. There are eight factors in the experiment: pulse rate ($A$), weld time ($B$), cool time ($C$), hold time ($D$), squeeze time ($E$), air pressure ($F$), current percentage ($G$), tip size ($H$). Among them, the pulse rate and the weld time are related factors, i.e., for lower pulse rate, the adequate weld time should be set longer in order to produce weld points with acceptable quality. An 18-run orthogonal array, $OA(18, 2^1 3^7)$, with a slight modification was adopted to study the eight factors. The planning matrix of these factors are given in Table 1 (unfortunately, the units of these factors was not reported). Factors $A$ and $H$ have two levels and other factors have three levels. Note that the column $H$ in Table 1 is obtained by collapsing a three-level factor in the $OA(18, 2^1 3^7)$ to a two-level factor (see Wu and Hamada, [8], Section 7.8). The actual settings of low, median, and high levels in the column $B$ depend on the levels of $A$ as shown in Table 2. We regard the area enclosed by solid lines in Fig 2 as the adequate experimental region, i.e., $R_E$.

TABLE 1
*Planning matrix of the welding experiment*

| A | B | C | D | E | F | G | H |
|---|---|---|---|---|---|---|---|
| 2 | low | 6 | 10 | 15 | 50 | 85 | 3/8 |
| 2 | low | 12 | 18 | 20 | 55 | 90 | 1/4 |
| 2 | low | 18 | 26 | 25 | 60 | 95 | 3/8 |
| 2 | median | 6 | 10 | 20 | 55 | 95 | 3/8 |
| 2 | median | 12 | 18 | 25 | 60 | 85 | 3/8 |
| 2 | median | 18 | 26 | 15 | 50 | 90 | 1/4 |
| 2 | high | 6 | 18 | 15 | 60 | 90 | 3/8 |
| 2 | high | 12 | 26 | 20 | 50 | 95 | 3/8 |
| 2 | high | 18 | 10 | 25 | 55 | 85 | 1/4 |
| 4 | low | 6 | 26 | 25 | 55 | 90 | 3/8 |
| 4 | low | 12 | 10 | 15 | 60 | 95 | 1/4 |
| 4 | low | 18 | 18 | 20 | 50 | 85 | 3/8 |
| 4 | median | 6 | 18 | 25 | 50 | 95 | 1/4 |
| 4 | median | 12 | 26 | 15 | 55 | 85 | 3/8 |
| 4 | median | 18 | 10 | 20 | 60 | 90 | 3/8 |
| 4 | high | 6 | 26 | 20 | 60 | 85 | 1/4 |
| 4 | high | 12 | 10 | 25 | 50 | 90 | 3/8 |
| 4 | high | 18 | 18 | 15 | 55 | 95 | 3/8 |

TABLE 2
*Actual settings of pulse rate and weld time*

| pulse rate | weld time | | |
|---|---|---|---|
| | low | median | high |
| 2 | 32 | 36 | 40 |
| 4 | 18 | 22 | 26 |

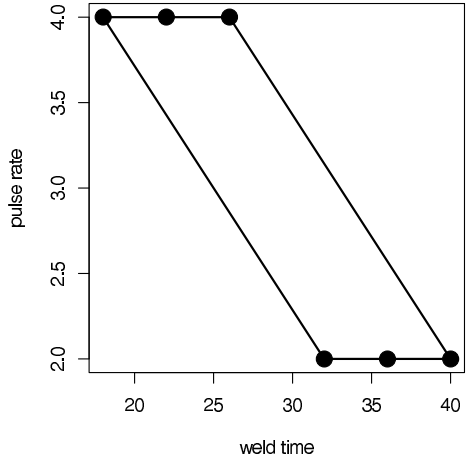

FIG 2. *Adequate experimental region of pulse rate and weld time.*

The main difference between the three modeling techniques is reflected in the effect coding of the slid factor $B$. For the RCRS model, because the $R_E$ is transformed into a square after re-centering and re-scaling, the effect coding of $B$ is the same as in a non-slid factor. Therefore, by applying the linear-quadratic system in [8], Section 5.6, the linear effect of $B$ codes the low, median, and high levels as $-1$, 0, and 1, respectively, and the quadratic effect of $B$ as 1, $-2$, 1, respectively. They are shown in the columns labeled by $B_l$ and $B_q$ of Table 3. Note that, although we still call $B_l$ and $B_q$ the main effects of $B$, they are no longer the main effects of weld time. Instead, the $B$ after RCRS represents a new factor which is a linear combination of weld time and pulse rate. For example, from $B_l = -1$ in Table 3, we can see that the low level of the new factor is the left hand side boundary of $R_E$ in Fig 2 (i.e., the straight line that links the point (weld time, pulse rate)=(32, 2) and the point (weld time, pulse rate)=(18, 4)) and from $B_l = 1$ the high level is the right hand side boundary. For the NEM approach, the effects of $B$ are conditional on the levels of $A$. For each level of $A$, the linear-quadratic system is applied to generate the $B_l|A_1$, $B_q|A_1$, $B_l|A_2$, and $B_q|A_2$ as shown in Table 3. For the RSM approach, because the lowest actual setting of $B$ is 18 and the highest actual setting is 40, we code 18 as $-1$ and 40 as $+1$, and the other settings, 22, 26, 32, and 36, are proportionally coded as $-\frac{7}{11}$, $-\frac{3}{11}$, $\frac{3}{11}$, and $\frac{7}{11}$, respectively. These are shown in the column labeled as $x_B$ in Table 3. The $x_B^2$ is the componentwise square of $x_B$.

TABLE 3
*Effect coding of pulse rate and weld time for the three modeling techniques*

| factors | | RCRS | | | NEM | | | | | RSM | | |
|---|---|---|---|---|---|---|---|---|---|---|---|---|
| $A$ | $B$ | $A_l$ | $B_l$ | $B_q$ | $A_l$ | $B_l\|A_1$ | $B_q\|A_1$ | $B_l\|A_2$ | $B_q\|A_2$ | $x_A$ | $x_B$ | $x_B^2$ |
| 2 | low(32)    | $-1$ | $-1$ | 1    | $-1$ | $-1$ | 1    | 0    | 0    | $-1$ | 0.273    | 0.074 |
| 2 | median(36) | $-1$ | 0    | $-2$ | $-1$ | 0    | $-2$ | 0    | 0    | $-1$ | 0.636    | 0.405 |
| 2 | high(40)   | $-1$ | 1    | 1    | $-1$ | 1    | 1    | 0    | 0    | $-1$ | 1        | 1 |
| 4 | low(18)    | 1    | $-1$ | 1    | 1    | 0    | 0    | $-1$ | 1    | 1    | $-1$     | 1 |
| 4 | median(22) | 1    | 0    | $-2$ | 1    | 0    | 0    | 0    | $-2$ | 1    | $-0.636$ | 0.405 |
| 4 | high(26)   | 1    | 1    | 1    | 1    | 0    | 0    | 1    | 1    | 1    | $-0.273$ | 0.074 |

In the data analysis, we consider the models that contain all main effects of factors $C$-$H$ and five effects generated from factors $A$ and $B$. For the RCRS, the five effects are $A_l$, $B_l$, $B_q$, $A_lB_l$, and $A_lB_q$, where $A_lB_l$ and $A_lB_q$ are interactions generated by the componentwise multiplication of $A_l$ and $B_l$, and $A_l$ and $B_q$, respectively. For the NEM, the five effects are $A_l$, $B_l|A_1$, $B_q|A_1$, $B_l|A_2$, and $B_q|A_2$. For the RSM, the five effects are $x_A$, $x_B$, $x_B^2$, $x_Ax_B$, and $x_Ax_B^2$, where $x_Ax_B$ and $x_Ax_B^2$ are interactions generated by the componentwise multiplication of $x_A$ and $x_B$, and $x_A$ and $x_B^2$, respectively. Although the five effects are coded in different ways for each modeling technique, the vector spaces spanned by any set of the five effects are identical. Consequently, the effects of factors $C$-$H$ will have the same analysis results in the three models. Because of this reason, we only give the analysis results of the five effects generated by $A$ and $B$, which include their estimated values, $t$-values, and $p$-values, under the RCRS, the NEM, and the RSM, in Tables 4, 5, and 6, respectively. From these tables, we have some interesting findings presented in the following.

1. The $B_l|A_1$ and $B_l|A_2$ are the linear effects of $B$ conditional on two different levels of $A$. In Table 5, we find that the two conditional effects have different magnitudes. When the pulse rate is 2, the weld time has a strong linear effect (significant $B_l|A_1$). When the pulse rate is changed to 4, the linear effect of weld time ($B_l|A_2$) is insignificant. After re-centering and re-scaling, these two effects are transformed into two parameters, $B_l$ and $A_lB_l$, in Table 4. The $B_l$ represents the average of the two conditional linear effects ($78.96 = ((-23.75) + 181.67)/2$) and the interaction $A_lB_l$ represents the difference between the two conditional linear effects ($-102.71 = ((-23.75) - 181.67)/2$). It is then clear why $B_l$ and $A_lB_l$ are both significant. The same argument can be applied to $B_q|A_1$ and $B_q|A_2$ in Table 5 and $B_q$ and $A_lB_q$ in Table 4. Because $B_q|A_1$ and $B_q|A_2$ has rather similar magnitudes ($-27.92$ and $-41.67$), it explains why their difference (i.e., $A_lB_q$) is insignificant.

2. By comparing $x_A$ in Table 6 and $A_l$ in Tables 4 and 5, we find surprisingly that $A_l$ is significant while $x_A$ is insignificant even though $A_l$ and $x_A$ have the same coding in Table 3. By a further investigation of the correlations between the estimated effects (given in Table 7), it is seen that the insignificance of $x_A$ is caused by the severe collinearity between $x_A$ and $x_B$. It also results in the other three high correlations in Table 7 because other effects are also defined by $x_A$ and $x_B$. Note that in the planning matrix in Table 1, all effects in the models based on the RCRS and the NEM are mutually orthogonal. The appearance of severe collinearity will be further explained in Section 5.

TABLE 4
*Analysis based on RCRS*

|       | value   | $t$-value | $p$-value |
|-------|---------|-----------|-----------|
| $A_l$ | $-81.04$ | $-6.20$ | 0.00 |
| $B_l$ | 78.96 | 4.93 | 0.00 |
| $B_q$ | $-34.79$ | $-3.76$ | 0.00 |
| $A_l B_l$ | $-102.71$ | $-6.42$ | 0.00 |
| $A_l B_q$ | $-6.88$ | $-0.74$ | 0.46 |

TABLE 5
*Analysis based on NEM*

|       | value   | $t$-value | $p$-value |
|-------|---------|-----------|-----------|
| $A_l$ | $-81.04$ | $-6.20$ | 0.00 |
| $B_l|A_1$ | 181.67 | 8.03 | 0.00 |
| $B_q|A_1$ | $-27.92$ | $-2.14$ | 0.04 |
| $B_l|A_2$ | $-23.75$ | $-1.05$ | 0.30 |
| $B_q|A_2$ | $-41.67$ | $-3.19$ | 0.00 |

TABLE 6
*Analysis based on RSM*

|       | value   | $t$-value | $p$-value |
|-------|---------|-----------|-----------|
| $x_A$ | 7.72 | 0.11 | 0.92 |
| $x_B$ | 18.62 | 0.07 | 0.95 |
| $x_B^2$ | $-789.34$ | $-3.76$ | 0.00 |
| $x_A x_B$ | $-1287.06$ | $-4.76$ | 0.00 |
| $x_A x_B^2$ | $-155.98$ | $-0.74$ | 0.46 |

TABLE 7
*Correlation matrix of the estimated effects under the RSM model*

|       | $x_B$ | $x_B^2$ | $x_A x_B$ | $x_A x_B^2$ |
|-------|-------|---------|-----------|-------------|
| $x_A$ | 0.96 | 0.00 | 0.00 | 0.91 |
| $x_B$ |      | 0.00 | 0.00 | 0.99 |
| $x_B^2$ |    |      | 0.99 | 0.00 |
| $x_A x_B$ |  |      |      | 0.00 |

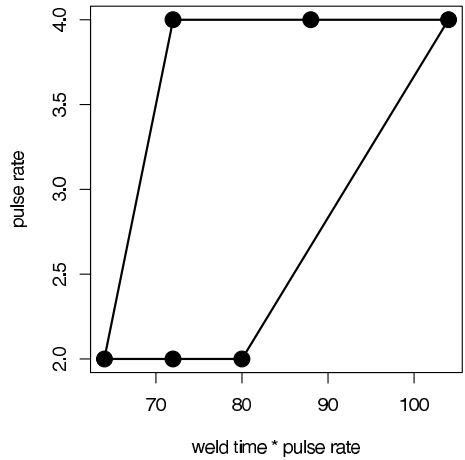

FIG 3. *Transformed experimental region.*

Suppose that severe collinearity is a serious concern but analysis based on RCRS or NEM is not an option to the investigators. A possible choice for reducing the collinearity might be to transform the variables. For example, after replacing weld time in Table 1 by a new variable, pulse rate times weld time, the RSM model will exhibit less correlation between the parameter estimates because the transformed experimental region (given in Fig 3) is more similar to a rectangle than the original experimental region (as shown in Fig 2).

3. From Table 3, we can understand that for the main effects of the slid factor (which include $B_l$ and $B_q$ in RCRS, $B_l|A_i$ and $B_q|A_i$ in NEM, and $x_B$ and $x_B^2$ in RSM), the coding based on RCRS can best preserve orthogonality property in a planning matrix, followed by the NEM, and RSM being the worst.

## 5. Relationship between RCRS and RSM models

To explain the relationship between the RCRS and RSM models, consider an RCRS model for the nine-run experiment that contains all main effects and a linear-by-linear interaction as follows:

$$
(7) \quad \begin{aligned} f(x_A, x_B) &\approx \eta_0 + \eta_1 x_A + \eta_{11} x_A^2 + \eta_2 \left[\frac{x_B - m_B(x_A)}{l_B(x_A)}\right] \\ &+ \eta_{22} \left[\frac{x_B - m_B(x_A)}{l_B(x_A)}\right]^2 + \eta_{12} x_A \left[\frac{x_B - m_B(x_A)}{l_B(x_A)}\right], \end{aligned}
$$

where $m_B(x_A)$ and $l_B(x_A)$ are the center and range, respectively, of the experimental region chosen for $B$ when $A$ is conditioned on $x_A$. For simplicity, assume that $m_B$ is a linear function of $x_A$, i.e., $m_B(x_A) = s + t\,x_A$, and $l_B$ is a constant, i.e., $l_B(x_A) = r$ (same shape as Fig 2). By substituting them into (7) and expanding (7) in a polynomial form, we obtain an RSM model, consisting of the factorial effects $x_A$, $x_A^2$, $x_B$, $x_B^2$, and $x_A x_B$, as follows:

$$
(8) \quad \begin{aligned} f(x_A, x_B) &\approx \left[\eta_0 + (s^2/r^2)\eta_{22} - (s/r)\eta_2\right] \\ &+ \left[\eta_1 - (t/r)\eta_2 - (s/r)\eta_{12} + (2st/r^2)\eta_{22}\right] x_A \\ &+ \left[\eta_{11} + (t^2/r^2)\eta_{22} - (t/r)\eta_{12}\right] x_A^2 + \left[(1/r)\eta_2 - (2s/r^2)\eta_{22}\right] x_B \\ &+ \left[(1/r^2)\eta_{22}\right] x_B^2 + \left[(1/r)\eta_{12} - (2t/r^2)\eta_{22}\right] x_A x_B. \end{aligned}
$$

Note that in (8), the parameters of factorial effects are functions of $\eta$'s and $r$, $s$, and $t$. The $\eta$'s represent the relationship between factors and response in the RCRS model, and $r$, $s$, and $t$ characterize the shape of the irregular experimental region. The shape has been eliminated in the RCRS model after applying the transformation $\frac{x_B - m_B(x_A)}{l_B(x_A)}$ on $B$. However, $m_B$ and $l_B$ still affect the polynomial terms of the RSM model in (8). This example shows that an RSM model for sliding-level experiments contains two components: a description of the relationship between factors and response, and a description of the irregular shape of the experimental region. The two components are intertwined and undistinguishable in the parameters of an RSM model. On the other hand, a fitted model based on RCRS only contains information on the first component because the irregularity of shape has been eliminated after re-centering and re-scaling. This observation is supported by the appearance of strong collinearity between $x_A$ and $x_B$ in Table 7. Note that after RCRS, the main effects of $A$ and $B$ (in Table 4) are orthogonal. However, in the RSM model such strong collinearity inevitably appears because: (i) the parameters in the RSM model are influenced by the irregular shape of experimental region, and (ii) the irregular shape (i.e., $R_E$ in Fig 2) reflects the fact that $B$ is smaller when $A$ is larger.

In general, the shape of the chosen experimental region can be arbitrary, and $m_B$ and $l_B$ can have more complicated forms than what was assumed above. However, similar remarks and conclusions are still applicable.

Suppose that the relationship between mean response and factors $A$ and $B$ satisfies (1). When the $m_B$ and $l_B$ in (7) are appropriately chosen so that the interaction elimination after RCRS is achieved, the $\eta_{12}$ in (7) becomes zero. In this case, the RCRS model in (7) does not nominally contain an interaction effect, but an interaction (i.e., $x_A x_B$) is still present in the RSM model (8). The apparent discrepancy lies in the different approaches they take to handle the irregular shape of the experimental region. This observation partially supports the interaction elimination rationale in RCRS from a different perspective. Because the fitted model after RCRS does not properly take into account the irregular shape of the experimental region, it can, in most cases, utilize fewer effects than an RSM model to achieve a comparable coefficient of determination (i.e., $R^2$). Some interaction effects (such as the $x_A x_B$ in the case) are not required for the RCRS model.

## 6. Summary

For the purpose of response prediction for sliding-level experiments, we point out the shortcomings of two existing approaches, RCRS and NEM, when the related factors are quantitative. An alternative analysis strategy is proposed based on the response surface modeling, in which the response prediction can be implemented in a straightforward manner. Through the comparisons of the three strategies, we present several interesting conclusions, which lead to better understanding of the concepts, properties, limitations, and implicit assumptions behind each strategy. None of the three methods dominates the others in every aspects. The best strategy for the investigators depends on the information they have about the irregular region and the objectives of the experiment. Although we do not discuss the design issues in this article, the choice of the modeling strategy will influence the choice of the best design. This and other issues in modeling deserve further study.

### Acknowledgments

S.-W. Cheng's research was supported by the National Science Council of Taiwan, ROC. C.F.J. Wu's research was supported by the NSF. L. Huwang's research was supported by the National Science Council of Taiwan, ROC, under the grant No. NSC-94-2118-M-007-003. The authors are grateful to the referee for valuable comments.